\documentclass{amsart}
\usepackage{amssymb}
\usepackage{amsfonts}

\theoremstyle{plain}

\numberwithin{equation}{section}
\input{tcilatex}

\begin{document}
\title[Balancedness]{Balancedness of Social Choice Correspondences}
\author{Jerry S. Kelly}
\address{JK: Department of Economics, Syracuse University, Syracuse NY
13244-1020 USA}
\email{jskelly@syr.edu}
\author{Shaofang Qi}
\address{SQ: Humboldt University Berlin, School of Business and Economics,
Spandauer Str. 1, 10178 Berlin, Germany}
\email{qi.shaofang@hu-berlin.de}
\date{March 9, 2018}

\begin{abstract}
A social choice correspondence satisfies \textbf{balancedness} if, for every
pair of alternatives, $x$ and $y$, and every pair of individuals, $i$ and $j$%
, whenever a profile has $x$ adjacent to but just above $y$ for individual $%
i $ while individual $j$ has $y$ adjacent to but just above $x$, then only
switching $x$ and $y$ in the orderings for both of those two individuals
leaves the choice set unchanged. We show how the balancedness condition
interacts with other social choice properties, especially tops-only. We also
use balancedness to characterize the Borda rule (for a fixed number of
voters) within the class of scoring rules.
\end{abstract}

\maketitle

\section{Introduction}

We consider the interaction of several social choice properties with a new
condition, balancedness. A social choice correspondence satisfies \textbf{%
balancedness} if, for every pair of alternatives, $x$ and $y$, and every
pair of individuals, $i$ and $j$, whenever a profile has $x$ adjacent to but
just above $y$ for individual $i$ while individual $j$ has $y$ adjacent to
but just above $x$, then only switching $x$ and $y$ in the orderings for
both of those two individuals leaves the choice set unchanged.\medskip

Social choice theory often considers responsiveness conditions, like
monotonicity, but balancedness is a \textit{non-responsiveness} property. \
It is a natural equity condition that simultaneously incorporates some equal
treatment for individuals short of anonymity, some equal treatment of
alternatives short of neutrality, and some equal treatment for differences
of position of alternatives in orderings (for example, raising $x$ just
above $y$ in the bottom two ranks for individual $j$ exactly offsets
lowering $x$ just below $y$ in the top two ranks for $i$).\medskip

Let $X$ with cardinality $|X|=m\geq 2$ be the set of \textbf{alternatives}
and let $N=\{1,2,...,n\}$ with $n\geq 2$ be the set of \textbf{individuals}.
A (strong) \textbf{ordering} on $X$ is a complete, asymmetric, transitive
relation on $X$ (non-trivial individual indifference is disallowed). The
highest ranked element of an ordering $r$ is denoted $r[1]$, the second
highest is denoted $r[2]$, etc. \ And\ $r[1:k]$ is the unordered set of
alternatives in the top $k$ ranks of $r$. \ The set of all orderings on $X$
is $L(X)$. A \textbf{profile} $u$ is an element $(u(1),u(2),...,u(n))$ of
the Cartesian product $L(X)^{N}$.\medskip 

A \textbf{social choice correspondence} $G$ is a map from the domain\ $%
L(X)^{N}$ to non-empty subsets of $X$. The \textbf{range} of social choice
correspondence $G$ is the collection of all sets $S$ such that there exists
a profile $u$ with $G(u)=S$.\medskip

$G$ satisfies \textbf{tops-only} if for all profiles $u$, $v$, whenever $%
u(i)[1]=v(i)[1]$ for all $i$, then $G(u)=G(v)$.\medskip

We say profile $v$ is \textbf{constructed from profile }$u$\textbf{\ by
transposition pair} $\mathbf{(x,y)}$\textbf{\ via individuals }$i$\textbf{\
and }$j$ if at $u$, $x$ is immediately above $y$ for $i$, and $y$ is
immediately above $x$ for $j$, and profile $v$ is just the same as $u$
except that for $i$ and $j$, alternatives $x$ and $y$ are transposed. A
social choice correspondence $G$ will be called \textbf{balanced} if, for
all $x$, $y$, $u$, $v$, $i$, and $j$, whenever profile $v$ is constructed
from $u$ by transposition pair $(x,y)$ via individuals $i$ and $j$, then $%
G(v)=G(u)$. (Otherwise, $G$ will be called \textbf{unbalanced)}.\medskip

If $m=2$, balancedness is equivalent to anonymity. For $m>2$, balancedness
holds for the Pareto correspondence\footnote{%
For definitions and discussions of these social choice rules and the
properties mentioned in this paragraph, see Arrow, Sen, and Suzumura (2002)
and Heckelman and Miller (2015).}, the Borda rule, the Copeland rule, and
top-cycle (which selects the maximal set of the transitive closure of simple
majority voting). (Thus balancedness is consistent with anonymity,
neutrality, and the Pareto condition.) Those rules all fail tops-only.\
Balancedness fails for the tops-only correspondences that are dictatorship,
plurality, and union-of-the-tops. (Thus unbalancedness is also consistent
with unanimity, anonymity, neutrality, and the Pareto condition.) \ In fact,
balancedness fails for almost all tops-only correspondences as Section 1
will show. Balancedness also fails for the maximin rule.\medskip

To see one way balancedness reflects equal treatment of individuals, first
observe again that any dictatorial correspondence (where for some $i$ and
all $u$, $G(u)=u(i)[1]$) is unbalanced. \ This can be extended to cover all
correspondences with ineffective individuals. \ An individual $i$ is \textbf{%
ineffective} for a correspondence $G$ if for all profiles $u$, $u^{\ast }$,
we have $G(u)=G(u^{\ast })$ whenever $u(j)=u^{\ast }(j)$ for all $j\neq i$.
\ For example, if $G$ is dictatorial with dictator $j$, then \textit{every} 
\textit{other} individual is ineffective. We show that balancedness implies
every individual is effective.\medskip

\textbf{Theorem 1.} \ Any non-constant social choice correspondence with an
ineffective individual is unbalanced.\medskip

\textbf{Proof:} \ Let $i$ be an ineffective individual for non-constant
correspondence $G$. \ Let $u$, $u^{\ast }$ be two profiles with $G(u)\neq
G(u^{\ast })$. \ Since $i$ is ineffective, we may assume that $u(i)=u^{\ast
}(i)$. \ Now construct a sequence of profiles

\begin{equation*}
u=u_{1},u_{2},...,u_{T-1},u_{T}=u^{\ast }
\end{equation*}%
\newline
from $u$ to $u^{\ast }$ such that for any two successive profiles $u_{t-1}$, 
$u_{t}$ in the sequence, $u_{t}$ differs from $u_{t-1}$ only by a
transposition of two adjacent alternatives in the ordering of a single
individual other than $i$. \ For some $t$, it must be that $G(u_{t})\neq
G(u_{t-1})$. \ Suppose that $u_{t}$ differs from $u_{t-1}$ because $x$ is
transposed with $y$, ranked just below $x$ in $u_{t-1}(j)$. \medskip

Construct profile $u_{t-1}^{\prime }$ from $u_{t-1}$ by moving $y$ adjacent
to and just above $x$ in $u_{t-1}(i)$, and construct profile $u_{t}^{\prime
} $ from $u_{t}$ by moving $x$ adjacent to and just above $y$ in $u_{t-1}(i)$%
. Then, by the ineffectiveness of $i$, $G$($u_{t}^{\prime })=G(u_{t})$ and $%
G(u_{t-1}^{\prime })=G(u_{t-1})$, so $G(u_{t}^{\prime })\neq
G(u_{t-1}^{\prime })$. \ But $u_{t}^{\prime }$ differs from $u_{t-1}^{\prime
}$ by transposition pair $(x,y)$ via individuals $i$ and $j$ and so this
constitutes a violation of balancedness. \ \ \ \ $\square $\medskip

In the first half of this paper, we show how the balancedness condition
interacts with other social choice properties, especially tops-only. \ In
the second half, we show how balancedness can be used to characterize the
Borda rule (for a fixed number of voters) within the class of scoring rules.

\section{Tops-only}

Balancedness incorporates some equal treatment of positions, so we might
expect conflict with the tops-only property. \ By tops-only, all profiles $u$
with an alternative $x$ at everyone's top must give the same value for $G(u)$%
, though it need not be the case that $G(u)=\{x\}$. \ Do there exist any
social choice correspondences satisfying tops-only and unanimity that are
balanced? \ A positive answer to that question is given by:\medskip

\textbf{Example 1.} \ Let $G(u)$ be the common top at all unanimous profiles
and set $G(u)=X$ (or any other fixed set) on all non-unanimous
profiles.\medskip

That example illustrates the following result.\medskip

\textbf{Theorem 2}. \ For $m\geq 3$, $n\geq 3$, if social choice
correspondence $G$ satisfies tops-only and balancedness, then it is constant
on all non-unanimous profiles.\medskip

Thus such $G$ has at most $m+1$ sets in its range.\medskip

\textbf{Proof:} \ Let $a$ and $b$ be two fixed elements of $X$ and let $%
u^{\ast }$ be a profile with tops $abb...bb$. \ Given any profile $u$
without a unanimous top, we will show there is a sequence of profiles $%
u_{1},u_{2},...,u_{T}$ all without a unanimous top, such that:\medskip

\qquad \qquad 1. \ $u_{1}=u$;

\qquad \qquad 2. \ $u_{T}=u^{\ast }$;

\qquad \qquad 3. \ For all $t$ with $1\leq t<T$, $u_{t+1}$ can be
constructed from $u_{t}$ by either a reordering of alternatives below
someone's top or by a paired transposition.\medskip

Then $G(u)=G(u^{\ast })$ for all $u$ without a unanimous top.\medskip

In the following, let $T(u)$ be the set of all alternatives $x$ such that $%
u(i)[1]=x$ for some $i$.\medskip

\textbf{Case 1}. \ $T(u)=\{a,b\}$. \ If $u(1)[1]=a$, go to the next
paragraph. \ Suppose that $u(1)[1]=b$. For some $i>1$ with $a$ in the top
rank, construct $u_{2}$ by raising $a$ to 1's second rank and $b$ to $i$'s
second rank. \ Then construct $u_{3}$ by transposition pair $(a,b)$ via
individuals 1 and $i$. \ Now $u_{3}(1)[1]=a$ and $u_{3}(i)[1]=b$.\medskip

If all other tops are now $b$, we are done. \ So suppose that some other
top, say for individual $j$, is $a$. \ Construct $u_{4}$ by raising a third
alternative $c$ to the second rank for $j$ and raising $c$ just above $a$
for $i$. \ Then construct $u_{5}$ by transposition pair $(a,c)$ via $i$ and $%
j$. Now $c$ is at $j$'s top. Recall that at $u_{5}$, we have $u_{5}(1)[1]=a$%
, so we can move $b$ to the second rank for individual 1, $c$ to the third
rank for individual 1, and for individual $j$, we move $b$ to the second
rank. Construct the next profile in the sequence by transposition pair $%
(b,c) $ via 1 and $j$. Then $b$ now becomes $j$'s top; the $a$ has been
changed to a $b$. Repeat this until all "$a$"s have been changed first to "$%
c $"s and then to "$b$"s.\ \medskip

\textbf{Case 2}. \ Suppose $T(u)$ contains $a$, $b$, and other alternatives.
\ If $c$ is a top for someone, say $i$, construct $u_{2}$ by raising $b$ to
the second rank for $i$ and raise $b$ just above $c$ for an individual $j$
who has $a$ on top. \ Then a transposition pair $(b,c)$ via $i$ and $j$
reduces by one the number of individuals without $a$ or $b$ on top. \
Continue in this fashion until you reach Case 1.\medskip

\textbf{Case 3}. \ $T(u)$ contains one of $a$ and $b$ but not the other. \
Suppose that $u(i)[1]=a$ and $u(j)[1]=c$. \ Construct $u_{2}$ by raising $b$
to second rank for $j$ and raise $b$ just above $c$ for $i$. \ Then a
transposition of $b$ and $c$ for $i$ and $j$ creates a profile in either
Case 2 or Case 1. \ A similar analysis holds if it is $b$ instead of $a$ at
someone's top.\medskip

\textbf{Case 4}. \ $T(u)$ contains neither of $a$ or $b$ but does contain
say $c$ and $d$. \ Construct profile $u_{2}$ by raising $a$ to second rank
for some $i$ with $c$ on top and $a$ just above $c$ for some $j$ with $d$ on
top. \ A transposition pair $(a,c)$ via $i$ and $j$ yields a profile in Case
3. \ \ \ \ $\square $\medskip

A modified version of this argument works also for $n=2$.\medskip

But there is a limit to this style of argument; there are not sequences to $%
u^{\ast }$ from profiles with unanimous tops as balancedness cannot be
applied there. This limit on the argument can not be overcome, as seen by
the example at the beginning of this section; there, at profiles with a
common top, different outcomes may occur.\medskip

\textbf{Corollary.} \ For $m\geq 3$, if social choice correspondence $G$
satisfies tops-only and Pareto, then $G$ is unbalanced.\medskip

Both tops-only and Pareto are needed in the Corollary. \ If Pareto is not
assumed, a constant correspondence satisfies both tops-only and
balancedness. \ If tops-only is not assumed, the correspondence that selects
the set of Pareto optimal alternatives satisfies Pareto and balancedness.

\section{Top-2\protect\medskip}

The following social choice correspondence violates tops-only, but outcomes
do depend only on which alternatives are ranked first or second by
individuals (but not on how the top alternatives are ordered within those
two ranks).\medskip

\textbf{Example 2.} \ (Rigid\footnote{%
"Rigid" because the number of ranks (here 2) is fixed to be the same for all
individuals. For approval voting without rigidity, see Brams and Fishburn
(1983). \ Rigid k-approval voting appeared in Alemante, Campbell, and Kelly
(2015) where it was called approval voting (type-k).}) $2$-approval voting
is a social choice correspondence $G$ that is like plurality rule except
that instead of selecting the alternatives with the most frequently
occurring tops, it selects the alternatives with the most frequent
occurrences in the top two ranks for everyone. \ Like plurality, it is
unbalanced.\medskip

We now extend Theorem 1 to cover social choice correspondences that, like
2-approval voting, depend only on the top two ranks for every individual. \
We first need to extend tops-only. \ $u(i)[1:2]$ is the (unordered) set of
alternatives in the top two ranks for individual $i$ at profile $u$. \ Then
we say social choice correspondence $G$ satisfies \textbf{top-2-only} if for
all profiles $u$, $u^{\ast }$, $G(u)=G(u^{\ast })$ if $u(i)[1:2]=u^{\ast
}(i)[1:2]$ for all individuals $i$.\medskip

Let \textbf{D} be the subdomain of $L(X)^{N}$ consisting of all profiles for
which it is\textbf{\ not} true that $u(i)[1:2]=u(j)[1:2]$ for all
individuals $i$ and $j$, that is, at least three alternatives occur in the
top two ranks over all individuals. \textbf{D }is the analog here of the
subdomain of non-unanimous profiles in Section 1.\medskip

\textbf{Theorem 3}. \ Let $n\geq 3$ and $m\geq 4$. \ Then any social choice
correspondence $G$ satisfying balancedness and top-2-only must be constant
on \textbf{D}.\medskip

\textbf{Proof:} \ Consider a specific profile $u^{\ast }$ with $c$ in
everyone's top rank, $a$ in \#1's second rank, and $b$ in everyone else's
second rank. \ It will suffice to show that for every profile $u$ in \textbf{%
D}, there is a sequence of profiles in \textbf{D} from $u$ to $u^{\ast }$
such that each is obtained from the previous profile by a transposition pair
or an application of top-2-only.\medskip

We argue by induction on the number $p$ of individuals who, at $u$, \textbf{%
do not} have $c$ in their top two ranks.\medskip

\textbf{Basis:} \ Suppose that we have $p=0$ at $u$, i.e., every individual
has $c$ in their top two ranks. \ Construct $u^{\prime }$ from $u$ by (only)
raising $c$ to everyone's top rank. \ Now consider the subdomain \textbf{D}$%
^{\prime }$ of \textbf{D} consisting of all profiles in \textbf{D }where
everyone has $c$ top-ranked. \ Note that for each such profile at least two
distinct alternatives must be second-ranked. \ Let $G^{\prime }$ be the
restriction of $G$ to \textbf{D}$^{\prime }$. \ Social choice correspondence 
$G^{\prime }$ on \textbf{D}$^{\prime }$ induces a correspondence on those
profiles on $X\backslash \{c\}$ which have non-unanimous tops. \ Then the
analysis in the tops-only section shows that there is a sequence of profiles
in \textbf{D}$^{\prime }$ from $u^{\prime }$ to $u^{\ast }$.\medskip

\textbf{Induction step:} \ Assume now that there is a non-empty set S of
alternatives such that for all profiles $v$ in \textbf{D} with the number of
individuals who do not have $c$ in their top two ranks being less than $p$
we have the same outcome: $G(v)=S$. Suppose that $u$ in \textbf{D}\ is a
profile where the number of individuals who at $u$ do not have $c$ in their
top two ranks is $p>0$. \ We will show there is a sequence of profiles (such
that each is obtained from the previous profile by a transposition pair or
an application of top-2-only) from $u$ to a profile $u^{\prime }$ in \textbf{%
D} with $p-1$ individuals who do not have $c$ in the top two ranks (so $%
G(u^{\prime })=S$). \ Without loss of generality suppose that it is the
first $q=n-p$ individuals with $c$ ranked in the top two ranks. \ We may
assume that $c$ has been raised to the top for each of those $q=n-p$
individuals so that $u$ is:\medskip

\qquad \qquad \qquad \qquad 
\begin{tabular}{|l|l|l|l|l|l|l|}
\hline
$1$ & $2$ & $\cdots $ & $q$ & $q+1$ & $\cdots $ & $n$ \\ \hline
$c$ & $c$ &  & $c$ & $x$ &  & $z$ \\ 
$\vdots $ & $\vdots $ &  & $\vdots $ & $y$ &  & $w$ \\ 
&  &  &  & $\vdots $ &  & $\vdots $ \\ \hline
\end{tabular}%
\medskip \newline
Here $x$, $y$, $z$ and $w$ are distinct from $c$ (though these four need not
all be distinct from one another).\medskip

\textbf{Case 1}. \ At least one of the alternatives, say $x$, in the top two
ranks for $q+1,...,n$ is ranked second by one of $1,...,q$. \ Without loss
of generality, let $q+1$ be an individual with $x$ in the top two and let 1
be the individual with $c$ on top and $x$ in the second rank:\medskip

\qquad \qquad \qquad \qquad 
\begin{tabular}{|l|l|l|l|l|l|l|}
\hline
$1$ & $2$ & $\cdots $ & $q$ & $q+1$ & $\cdots $ & $n$ \\ \hline
$c$ & $c$ &  & $c$ & $x$ &  & $z$ \\ 
$x$ & $\vdots $ &  & $\vdots $ & $y$ &  & $w$ \\ 
$\vdots $ &  &  &  & $\vdots $ &  & $\vdots $ \\ \hline
\end{tabular}%
\medskip \newline
Then if necessary, lower $x$ to the second rank for $q+1$ and raise $c$ to
the third rank for that individual.\medskip

\qquad \qquad \qquad \qquad 
\begin{tabular}{|l|l|l|l|l|l|l|}
\hline
$1$ & $2$ & $\cdots $ & $q$ & $q+1$ & $\cdots $ & $n$ \\ \hline
$c$ & $c$ &  & $c$ & $y$ &  & $z$ \\ 
$x$ & $\vdots $ &  & $\vdots $ & $x$ &  & $w$ \\ 
$\vdots $ &  &  &  & $c$ &  & $\vdots $ \\ \hline
\end{tabular}%
\medskip \newline
Now a transposition pair $(c,x)$ via individuals 1 and $q+1$ yields a
profile in \textbf{D }with $q+1$ individuals with $c$ in the top two ranks
and so only $p-1$ individuals with $c$ not in the top two ranks.\medskip

\textbf{Case 2.} \ None of the alternatives in the top two ranks for $q+1$%
,...,$n$ is ranked second by one of 1,...,$q$.\medskip

\textbf{Subcase 2A.} \ There are two individuals among $1,...,q$ with the
same alternative $t$ ranked second:\medskip

\qquad \qquad \qquad \qquad 
\begin{tabular}{|l|l|l|l|l|l|l|}
\hline
$1$ & $2$ & $\cdots $ & $q$ & $q+1$ & $\cdots $ & $n$ \\ \hline
$c$ & $c$ &  & $c$ & $x$ &  & $z$ \\ 
$t$ & $t$ &  & $\vdots $ & $y$ &  & $w$ \\ 
$\vdots $ & $\vdots $ &  &  & $\vdots $ &  & $\vdots $ \\ \hline
\end{tabular}%
\medskip \newline
For one of those individuals, say \#2, move $y$, an alternative in $%
u(q+1)[1:2]$ (moved for $q+1$ to second rank if necessary) up to the third
rank and then raise $t$ to the third rank for $q+1$:\medskip

\qquad \qquad \qquad \qquad 
\begin{tabular}{|l|l|l|l|l|l|l|}
\hline
$1$ & $2$ & $\cdots $ & $q$ & $q+1$ & $\cdots $ & $n$ \\ \hline
$c$ & $c$ &  & $c$ & $x$ &  & $z$ \\ 
$t$ & $t$ &  & $\vdots $ & $y$ &  & $w$ \\ 
$\vdots $ & $y$ &  &  & $t$ &  & $\vdots $ \\ \hline
\end{tabular}%
\medskip \newline
Now a transposition pair $(t,y)$ via individuals 2 and $q+1$ yields a
profile in \textbf{D }still with $p$ individuals with $c$ not in the top two
ranks but now back in Case 1 (since $t$ is in the top two ranks for
individual $q+1$ and is ranked second for individual \#1).\medskip

\textbf{Subcase 2B.} \ No element in the top two ranks for $q+1,...,n$ is
second for any of $1,...,q$, and no two of $1,...,q$ have the same second
element. (Recall $p=n-q$.)\medskip

\textbf{Subsubcase 2Bi.} \ $\mathbf{p<n-1}$\textbf{\ (and so}\ $\mathbf{q>1}$%
\textbf{).}\medskip

\qquad \qquad \qquad \qquad 
\begin{tabular}{|l|l|l|l|l|l|l|}
\hline
$1$ & $2$ & $\cdots $ & $q$ & $q+1$ & $\cdots $ & $n$ \\ \hline
$c$ & $c$ &  & $c$ & $x$ &  & $z$ \\ 
$s$ & $t$ &  & $\vdots $ & $y$ &  & $w$ \\ 
$\vdots $ & $\vdots $ &  &  & $\vdots $ &  & $\vdots $ \\ 
&  &  &  &  &  &  \\ \hline
\end{tabular}%
\medskip \newline
Construct:\medskip

\qquad \qquad \qquad \qquad 
\begin{tabular}{|l|l|l|l|l|l|l|}
\hline
$1$ & $2$ & $\cdots $ & $q$ & $q+1$ & $\cdots $ & $n$ \\ \hline
$c$ & $c$ &  & $c$ & $x$ &  & $z$ \\ 
$s$ & $t$ &  & $\vdots $ & $y$ &  & $w$ \\ 
$t$ & $\vdots $ &  &  & $t$ &  & $\vdots $ \\ 
$y$ &  &  &  & $\vdots $ &  &  \\ \hline
\end{tabular}%
\medskip \newline
where t is 2's second and y is in q+1's top two (moved to second rank if
necessary).\medskip

Now a transposition pair $(t,y)$ via 1 and $q+1$ takes us back to Case 1
(since $t$ will be in individual $q+1$'s top two and in \#2's second
rank).\medskip

\textbf{Subsubcase 2Bii. \ }$\mathbf{p=n-1}$\textbf{\ (so }$\mathbf{q=1}$%
\textbf{)}.\medskip

\qquad \qquad \qquad \qquad 
\begin{tabular}{|l|l|l|l|}
\hline
$1$ & $2$ & $\cdots $ & $n$ \\ \hline
$c$ & $x$ &  & $z$ \\ 
$s$ & $y$ &  & $w$ \\ 
$\vdots $ & $\vdots $ &  & $\vdots $ \\ 
&  &  &  \\ \hline
\end{tabular}%
\medskip \newline
If not all of $2,...,n$ have the same top two alternatives, say $n$ has $z$
in the top two but 2 does not, move $z$ to second place for individual $n$
and raise $c$ to $n$'s third rank while raising $c$ and $z$ to third and
fourth rank for \#2:\medskip

\qquad \qquad \qquad \qquad 
\begin{tabular}{|l|l|l|l|}
\hline
$1$ & $2$ & $\cdots $ & $n$ \\ \hline
$c$ & $x$ &  & $w$ \\ 
$s$ & $y$ &  & $z$ \\ 
$\vdots $ & $c$ &  & $c$ \\ 
& $z$ &  & $\vdots $ \\ \hline
\end{tabular}%
\medskip \newline
Now a transposition pair $(c,z)$ via 2 and $n$ takes us back to $p=n-2$ in 
\textbf{D}.\medskip

On the other hand if $2,...,n$ all have the same two alternatives in the top
two ranks:\medskip

\qquad \qquad \qquad \qquad 
\begin{tabular}{|l|l|l|l|}
\hline
$1$ & $2$ & $\cdots $ & $n$ \\ \hline
$c$ & $x$ &  & $x$ \\ 
$s$ & $y$ &  & $y$ \\ 
$\vdots $ & $\vdots $ &  & $\vdots $ \\ 
&  &  &  \\ \hline
\end{tabular}%
\medskip \newline
Raise y to third rank for \#1 and s to third rank for \#2\medskip

\qquad \qquad \qquad \qquad 
\begin{tabular}{|l|l|l|l|}
\hline
$1$ & $2$ & $\cdots $ & $n$ \\ \hline
$c$ & $x$ &  & $x$ \\ 
$s$ & $y$ &  & $y$ \\ 
$y$ & $s$ &  & $\vdots $ \\ 
$\vdots $ & $\vdots $ &  &  \\ \hline
\end{tabular}%
\medskip \newline
Now a transposition pair $(s,y)$ via individuals 1 and $2$ takes us back (in 
\textbf{D}) to Case 1.\medskip

\textbf{Case 3. \ }$\mathbf{p=n}$\textbf{\ (so }$\mathbf{q=0}$\textbf{)}%
.\medskip

Because we are in \textbf{D}, there must be an alternative $z$, distinct
from $x$ and $y$ (\#1's top two), in the top two ranks for someone, say \#2.
\ Move $z$ to second rank for \#2,\medskip

\qquad \qquad \qquad \qquad 
\begin{tabular}{|l|l|l|l|}
\hline
$1$ & $2$ & $\cdots $ & $n$ \\ \hline
$x$ &  &  &  \\ 
$y$ & $z$ &  &  \\ 
$\vdots $ & $\vdots $ &  &  \\ 
&  &  &  \\ \hline
\end{tabular}%
\medskip \newline
Then raise $c$ and $z$ to ranks three and four for \#1 and $c$ to third rank
for \#2:\medskip

\qquad \qquad \qquad \qquad 
\begin{tabular}{|l|l|l|l|}
\hline
$1$ & $2$ & $\cdots $ & $n$ \\ \hline
$x$ &  &  &  \\ 
$y$ & $z$ &  &  \\ 
$c$ & $c$ &  &  \\ 
$z$ & $\vdots $ &  &  \\ \hline
\end{tabular}%
\medskip \newline
Now transposition pair $(c,z)$ via 1 and $2$ takes us back (in \textbf{D})
to Case 2Bii. \ \ \ $\square $\medskip \medskip

There are straightforward generalizations to top-3-only, etc.

\section{Balancedness and Borda}

Balancedness is an equity condition that incorporates some equal treatment
for differences of position of alternatives in orderings. \ As observed at
the beginning of this paper, raising $x$ just above $y$ in the bottom two
ranks for individual $j$ exactly offsets lowering $x$ just below $y$ in the
top two ranks for $i$. \ This equal treatment of differences of position
suggests trying to characterize the Borda rule within the class of scoring
rules. \ Let a \textbf{scoring system} be given by weights\medskip

\begin{equation*}
s_{1}\leq s_{2}\leq s_{3}\leq \cdots \leq s_{m}
\end{equation*}%
\medskip \newline
At profile $u$, the score for an alternative $x$ is the sum

\begin{equation*}
S(x,u)=\dsum\limits_{i=1}^{n}s(u,i,x)
\end{equation*}%
where $s(u,i,x)=s_{k}$ if $u(i)[k]=x$, i.e., the score for $x$ is the sum of
the weights corresponding to the ranks that $x$ occupies in the individual
orderings at $u$. \ The related \textbf{scoring social choice correspondence}
$G$ selects at $u$ the alternatives $x$ with lowest $S(x,u)$ values. \ If

\begin{equation*}
s_{1}=s_{2}=s_{3}=\cdots =s_{m}
\end{equation*}%
then $G$ is the constant social choice correspondence with $G(u)=X$ at all $%
u $, which is not very helpful. \ Accordingly, we henceforth only consider
systems of weights such that at least two weights are distinct.\medskip

The \textbf{Borda correspondence}, $G_{B}$, uses weights $1<2<3<\cdots <m$.
\ But of course other weights also generate Borda. \ If social choice
correspondence $G$ is generated by weights\medskip

\begin{equation*}
s_{1}\leq s_{2}\leq s_{3}\leq \cdots \leq s_{m}
\end{equation*}%
\newline
then $G$ is also generated by linearly transformed weights\medskip

\begin{equation*}
t_{1}\leq t_{2}\leq t_{3}\leq \cdots \leq t_{m}
\end{equation*}%
\medskip \newline
where $t_{i}=\alpha +\beta s_{i}$ for real numbers $\alpha $, $\beta $, with 
$\beta >0$.\medskip

\textbf{Lemma 0.} \ For $m\geq 3$ and $n\geq 2$: if $G$ is a scoring social
choice correspondence for scoring system given by\medskip

\begin{equation*}
s_{1}\leq s_{2}\leq s_{3}\leq \cdots \leq s_{m}
\end{equation*}%
\newline
with at least two weights distinct and if $G$ satisfies balancedness, then $%
s_{1}\neq s_{2}$.\medskip

\textbf{Proof}: \ If $s_{1}=s_{2}$, then there exists a $k$, $2\leq k<m$,
such that $s_{1}=s_{2}=...=s_{k}<s_{k+1}$. \ Construct profile $u$ such that 
$u(i)[1]=y$ and $u(i)[2]=x$ for all $i<n$ while $u(n)[k]=x$ and $u(n)[k+1]=y$%
. \ Then $x\in G(u)$ and $y\notin G(u)$. \ If profile $u^{\ast }$ is
constructed from $u$ by transposition pair $(x,y)$ via individuals $1$ and $%
n $, then $x\notin G(u^{\ast })$ and $y\in G(u^{\ast })$. \ So $G(u^{\ast
})\neq G(u)$ and $G$ fails balancedness. \ \ \ $\square $\medskip

As a consequence of Lemma 0, we assume from now on that the weights have
been transformed so that $s_{1}=1$ and $s_{2}=2$.\medskip

For profiles of strong orderings, we want to show that generally if a
scoring social choice correspondence satisfies balancedness, it must be the
Borda rule. \ However, the next section shows a limitation on this
objective.\medskip

\section{Borda: m = n = 3\protect\medskip}

\textbf{Example 3}: \ Let $m=n=3$ and set scoring weights to be $1,2$, and $%
3.1$.\medskip

The related correspondence, $G$, differs from Borda. \ At profile $u$%
:\medskip

\qquad \qquad \qquad \qquad 
\begin{tabular}{|l|l|l|}
\hline
$1$ & $2$ & $3$ \\ \hline
$x$ & $x$ & $y$ \\ 
$y$ & $y$ & $z$ \\ 
$z$ & $z$ & $x$ \\ \hline
\end{tabular}%
\medskip \newline
the Borda rule has a tie between $x$ and $y$ and $z$ is Pareto-dominated by $%
y$ so\ $G_{B}(u)=\{x,y\}$. \ But with scoring weights $1,2$, and $3.1$, the
scores of $x$ and $y$ are 5.1 and 5 respectively: $G(u)=\{y\}$. \ $G\neq
G_{B}$.\medskip

Nevertheless, $G$ is balanced, as can easily be checked. \ So, for $m=n=3$,
balancedness of a scoring rule does \textit{not} imply Borda.\medskip

This use of small numbers of individuals and alternatives is critical in
Example 3, as we will see.\medskip

\section{Borda: n \TEXTsymbol{>} 3\protect\medskip}

Our analysis proceeds by induction on m. \ We begin by looking at $m=3$ and
a few small values of $n>3$.\medskip

\textbf{Lemma 1}. \ For $m=3$ and $n=4$, 5, or 6: if $G$ is a scoring social
choice correspondence and $G$ satisfies balancedness, then $G$ is the Borda
correspondence.\medskip

\textit{Proof}: We show the following: if $G$ is a scoring rule but not the
Borda correspondence, then $G$ violates balancedness.\medskip

\textbf{Basis case, }$m=3$,\textbf{\ for Lemma 1}\medskip

($n=5$) \ We first treat the case $n=5$. \ Consider the following profile $v$%
:\medskip

\qquad \qquad \qquad \qquad 
\begin{tabular}{|l|l|l|l|l|}
\hline
$1$ & $2$ & $3$ & $4$ & $5$ \\ \hline
$x$ & $z$ & $y$ & $x$ & $z$ \\ 
$y$ & $y$ & $x$ & $z$ & $y$ \\ 
$z$ & $x$ & $z$ & $y$ & $x$ \\ \hline
\end{tabular}%
\medskip \newline

Consider a scoring rule $G$ with $1<2\leq s_{3}$ (using Lemma 0). \ We
examine the scores for $x,y,z$:

For $x$, $4+2s_{3}$;

For $z$, $4+2s_{3}$;

For $y$, $7+s_{3}$.\medskip

So under this scoring rule, either $x$ and $y$ are both in $G(v)$ or neither 
$x$ nor $y$ are in $G(v)$.\medskip

Now we consider a transposition pair $(x,y)$ via individuals 1 and 2 that
yields the profile $u$:\medskip

\qquad \qquad \qquad \qquad 
\begin{tabular}{|l|l|l|l|l|}
\hline
$1$ & $2$ & $3$ & $4$ & $5$ \\ \hline
$y$ & $z$ & $y$ & $x$ & $z$ \\ 
$x$ & $x$ & $x$ & $z$ & $y$ \\ 
$z$ & $y$ & $z$ & $y$ & $x$ \\ \hline
\end{tabular}%
\medskip \newline
We examine again the scores for $x,y,z$:

For $x$, $7+s_{3}$;

For $y$, $4+2s_{3}$;

For $z$, $4+2s_{3}$.\medskip

This time, either $y$ and $z$ are both in $G(u)$ or both are out of $G(u)$.
\ But by balancedness, $G(u)=G(v)$, so at each profile all scores must be
the same. \ Therefore $7+s_{3}=4+2s_{3}$, which has the unique solution $%
s_{3}=3$, i.e., $G$ is the Borda correspondence.\medskip

($n=4$) \ Next, we consider $n=4$. We examine the following profile $v$%
.\medskip

\qquad \qquad \qquad \qquad 
\begin{tabular}{|l|l|l|l|}
\hline
$1$ & $2$ & $3$ & $4$ \\ \hline
$x$ & $z$ & $x$ & $z$ \\ 
$y$ & $y$ & $y$ & $y$ \\ 
$z$ & $x$ & $z$ & $x$ \\ \hline
\end{tabular}%
\medskip \newline
Consider a scoring rule $G$ with weights: $1<2\leq s_{3}$ (using Lemma 0).
We examine the scores for $x,y,z$:

For $x$, $2+2s_{3}$;

For $z$, $2+2s_{3}$;

For $y$, $8$.\medskip

Under this scoring rule, either $x$ and $z$ are both in $G(v)$ or neither $x$
nor $z$ is in $G(v)$.\medskip

Now we consider a \textit{sequence} of transposition pairs $(x,y)$ first via
individuals 1 and 2, and then via 3 and 4, which yields the profile $u$%
:\medskip

\qquad \qquad \qquad \qquad 
\begin{tabular}{|l|l|l|l|}
\hline
$1$ & $2$ & $3$ & $4$ \\ \hline
$y$ & $z$ & $y$ & $z$ \\ 
$x$ & $x$ & $x$ & $x$ \\ 
$z$ & $y$ & $z$ & $y$ \\ \hline
\end{tabular}%
\medskip \newline
We examine again the scores for $x,y,z$:

For $x$, $8$;

For $y$, $2+2s_{3}$;

For $z$, $2+2s_{3}$.\medskip

This time, either $y$ and $z$ are both in $G(v)$ or neither $y$ nor $z$ is
in $G(v)$. \ But $u$ is obtained from $v$ by a sequence of transposition
pairs, and balancedness implies $G(u)=G(v)$, which means both must be $%
\{x,y,z\}$, which, in turn, means at each profile, all three scores must be
the same. \ Therefore, $2+2s_{3}=8$, i.e., $G$ is the Borda
correspondence.\medskip

($n=6$) \ Finally, the analysis for $n=6$ proceeds just as for $n=4$, but
starting from the profile $v$:\medskip

\qquad \qquad \qquad \qquad 
\begin{tabular}{|l|l|l|l|l|l|}
\hline
$1$ & $2$ & $3$ & $4$ & $5$ & $6$ \\ \hline
$y$ & $z$ & $y$ & $z$ & $y$ & $z$ \\ 
$x$ & $x$ & $x$ & $x$ & $x$ & $x$ \\ 
$z$ & $y$ & $z$ & $y$ & $z$ & $y$ \\ \hline
\end{tabular}%
\medskip \newline
and this time doing three transposition pairs (for 1 and 2, then 3 and 4,
and then 5 and 6).\medskip

For $n>6$, merely observe that any such $n$ is equal to one of 4, 5, or 6
plus some multiple of 3. \ Using profiles above that are expanded by that
many multiples of a voting paradox profile (on which all alternatives have
the same score regardless of the weighting scheme) provides profiles showing 
$s_{3}=3$. \ \ \ \ $\square $\medskip

\textbf{Lemma 2}. \ For $m\geq 3$ and $n=4$, 5, or 6: if $G$ is a scoring
social choice correspondence and $G$ satisfies balancedness, then $G$ is the
Borda correspondence.\medskip

\textit{Proof}: \ We prove by induction on $m$. \ The basis step is given by
Lemma 1. \ For $M\geq 4$, suppose that the lemma holds for all $m<M$ and we
have linearly transformed weights\medskip

\qquad \qquad $1<2\leq s_{3}\leq \cdots \leq s_{M}$\medskip \newline
such that for some $j>3$, it is the case that $s_{j}\neq j$. \ Let $G$ be
the corresponding scoring social choice correspondence. \ We show that there
is a profile $u$ such that $G(u)\neq G_{B}(u)$ and that $G$ must fail
balancedness. \ Consider the smallest integer $j$ such that $s_{j}\neq j$. \
If $j<M$, just take a profile $v$ that works for $j$ alternatives (using the
induction hypothesis) and append $M-j$ additional alternatives to everyone's
bottom. \ \medskip

So we need only consider the case where the first $j$ such that $s_{j}\neq j$
is $j=M$. \ The weights are:\medskip

\qquad \qquad $1<2<3<\cdots <M-1\leq w$ where $w\neq M$.\medskip \newline

Analysis here is for $w>M$. The same profile applies to the case $M-1\leq
w<M $. We treat $n=5$, but the same construction works for 4 or 6
individuals. (See the footnote below, for example, for $n=4$.) \ For $n=5$
we previously looked at profile:\medskip

\qquad \qquad \qquad \qquad 
\begin{tabular}{|l|l|l|l|l|}
\hline
$1$ & $2$ & $3$ & $4$ & $5$ \\ \hline
$x$ & $z$ & $y$ & $x$ & $z$ \\ 
$y$ & $y$ & $x$ & $z$ & $y$ \\ 
$z$ & $x$ & $z$ & $y$ & $x$ \\ \hline
\end{tabular}%
\medskip \newline
Now at each stage of the construction, we add another alternative, placing
it just above $y$ and $x$ for \#2, just below $x$ and $y$ for \#1, and at
the bottom for everyone else. At the first stage we insert $a$:\medskip

\qquad \qquad \qquad \qquad 
\begin{tabular}{|l|l|l|l|l|}
\hline
$1$ & $2$ & $3$ & $4$ & $5$ \\ \hline
$x$ & $z$ & $y$ & $x$ & $z$ \\ 
$y$ & $a$ & $x$ & $z$ & $y$ \\ 
$a$ & $y$ & $z$ & $y$ & $x$ \\ 
$z$ & $x$ & $a$ & $a$ & $a$ \\ \hline
\end{tabular}%
\medskip \newline
(after which \#1 and \#2 have opposite rankings).\footnote{%
So for $n=4$, we also look at the profile for $n=4$ in the basis step in
Lemma 1, then insert $a$ just above $y$ and $x$ for \#2, just below $x$ and $%
y$ for \#1, and at the bottom for \#3 and \#4.} \ The Borda scores are 11
for $x$, $y$, and $z$ and 17 for $a$: $G_{B}(u)=\{x,y,z\}$. \ For the
correspondence $G$ with weights 1, 2, 3, and $w_{4}>4$, the scores are $%
7+w_{4}$ for $x$ and $z$, 11 for $y$ and $5+3w_{4}$ for $a$. \ Because $%
w_{4}>4$, the smallest of these is $11$ and $G(u)=\{y\}$. \ Construct $v$ by
transposing $x$ and $y$ for \#1 and \#2; this yields score 11 for $x$, and $%
7+w_{4}$ for both $y$ and $z$ and $5+3w_{4}$ for $a$. \ Again the smallest
of these is $11$, so $G(v)=\{x\}$, a failure of balancedness.\medskip

Eventually the profile looks like $v^{\ast }$:\medskip

\qquad \qquad \qquad \qquad 
\begin{tabular}{|l|l|l|l|l|}
\hline
$1$ & $2$ & $3$ & $4$ & $5$ \\ \hline
$x$ & $z$ & $y$ & $x$ & $z$ \\ 
$y$ & $a_{M-3}$ & $x$ & $z$ & $y$ \\ 
$a_{1}$ & $\vdots $ & $z$ & $y$ & $x$ \\ 
$a_{2}$ & $a_{2}$ & $a_{1}$ & $a_{1}$ & $a_{1}$ \\ 
$\vdots $ & $a_{1}$ & $a_{2}$ & $a_{2}$ & $a_{2}$ \\ 
$a_{M-3}$ & $y$ & $\vdots $ & $\vdots $ & $\vdots $ \\ 
$z$ & $x$ & $a_{M-3}$ & $a_{M-3}$ & $a_{M-3}$ \\ \hline
\end{tabular}%
\medskip \newline
(after which \#1 and \#2 have opposite rankings). \ It is easy to check
that: $G_{B}(v^{\ast })=\{x,y,z\}$. \ For the correspondence $G$ with
weights 1, 2, 3, ..., $w>M$, the scores are $7+w$ for $x$ and $z$, $7+M$ for 
$y$, and the rest $a_{1},a_{2},...$ have higher scores. \ Because $w>M$, the
smallest of these is $7+M$ and $G(v^{\ast })=\{y\}$.\medskip

Construct $u^{\ast }$ by transposing $x$ and $y$ for \#1 and \#2. \ This
yields scores $7+M$ for $x$, and $7+w$ for both $y$ and $z$, with the
remaining scores unchanged. \ Since $w>M$, we now have $G(u^{\ast })=\{x\}$,
a failure of balancedness. \ \ \ \ $\square $\medskip

\section{Borda: n = 3, m \TEXTsymbol{>} 3\protect\medskip}

Let's return to the case $n=3$, where we learned in Example 1 that for $m=3$
not every scoring social choice correspondence satisfying balancedness is
the Borda rule. \ For the case $m=4$ and $n=3$, a natural analog of Example
1 has scoring weights 1, 2, 3, and $4.1$. \ But that rule is unbalanced, as
can be seen at profile $u$:\medskip

\qquad \qquad \qquad \qquad 
\begin{tabular}{|l|l|l|}
\hline
$1$ & $2$ & $3$ \\ \hline
$a$ & $b$ & $x$ \\ 
$y$ & $a$ & $y$ \\ 
$x$ & $x$ & $b$ \\ 
$b$ & $y$ & $a$ \\ \hline
\end{tabular}%
\medskip \newline
Here the scores are $S(a)=7.1$, $S(b)=8.1$, $S(x)=7$, $S(y)=8.1$. \ The
smallest of these scores is 7 and so $G(u)=\{x\}$. \ If $v$ is constructed
from profile\textbf{\ }$u$\textbf{\ }by transposition pair $(x,y)$\ via
individuals $1$\ and $2$, then the scores become $S(a)=7.1$, $S(b)=8.1$, $%
S(x)=7.1$, $S(y)=8$. \ So $G(v)=\{x,a\}$, a failure of balancedness.\medskip

\textbf{Theorem 5}. \ For $n=3$ and $m>3$: if $G$ is a scoring social choice
correspondence and $G$ satisfies balancedness, then $G$ is the Borda
correspondence.\medskip

\textit{Proof}: \ As in the proof of Lemma 2, we argue by induction on $m$.
\ Starting with a basis case, new alternatives are inserted into a profile
to unbalancedness of a scoring rule with weights different from those of the
Borda rule. \ For the basis case here, with $m=4$, the scoring rule has
weights 1, 2, $s_{3}$, $s_{4}$.\medskip

\textbf{Basis case, }$m=4$,\textbf{\ for Theorem 5}\medskip

\textit{Proof.} Let's consider a scoring rule $G$ but not the Borda
correspondence: $1<2\leq s_{3}\leq s_{4}$. We consider the following profile 
$u$:\medskip

\qquad \qquad \qquad \qquad 
\begin{tabular}{|l|l|l|}
\hline
$1$ & $2$ & $3$ \\ \hline
$a$ & $b$ & $x$ \\ 
$y$ & $a$ & $y$ \\ 
$x$ & $x$ & $b$ \\ 
$b$ & $y$ & $a$ \\ \hline
\end{tabular}%
\medskip \newline
We examine the scores for $a,b,x,y$:\medskip

For $a$, $3+s_{4}$;

For $x$, $1+2s_{3}$;

For $b$, $1+s_{3}+s_{4}$;

For $y$, $4+s_{4}$.\medskip

Since $S(a)<S(y)$ and $S(x)<S(b)$ under this scoring rule, at least $a\in
G(u)$ or $x\in G(u)$ (or both).\medskip

Now we consider a transposition pair $(x,y)$ via individual 1 and 2 that
yields the profile $v$:\medskip

\qquad \qquad \qquad \qquad 
\begin{tabular}{|l|l|l|}
\hline
$1$ & $2$ & $3$ \\ \hline
$a$ & $b$ & $x$ \\ 
$x$ & $a$ & $y$ \\ 
$y$ & $y$ & $b$ \\ 
$b$ & $x$ & $a$ \\ \hline
\end{tabular}%
\medskip \newline
We examine the scores for $a,b,x,y$:\medskip

For $a$, $3+s_{4}$;

For $x$, $3+s_{4}$;

For $b$, $1+s_{3}+s_{4}$;

For $y$, $2+2s_{3}$.\medskip

By balancedness, $G(u)=G(v)$. Since at profile $v$, alternatives $a$ and $x$
have the same score, and at $u$, at least $a\in G(u)$ or $x\in G(u)$, we
have $a$ and $x$ are both in $G(u)$ and $G(v)$. So (from $u$), $%
S(a)=3+s_{4}=S(x)=1+2s_{3}$, i.e.,

$\qquad \qquad \qquad \qquad \qquad 2s_{3}-s_{4}=2$ \ \ \ \ (1)\medskip

Again, consider a transposition pair $(x,y)$ via individual 2 and 3 that
yields the profile $v^{\prime }$:\medskip

\qquad \qquad \qquad \qquad 
\begin{tabular}{|l|l|l|}
\hline
$1$ & $2$ & $3$ \\ \hline
$a$ & $b$ & $y$ \\ 
$x$ & $a$ & $x$ \\ 
$y$ & $x$ & $b$ \\ 
$b$ & $y$ & $a$ \\ \hline
\end{tabular}%
\medskip \newline
We examine the scores for $a,b,x,y$:\medskip

For $a$, $3+s_{4}$;

For $x$, $4+s_{3}$;

For $b$, $1+s_{3}+s_{4}$;

For $y$, $1+s_{3}+s_{4}$\medskip

By balancedness, $a,x\in G(v^{\prime })$, so $4+s_{3}=3+s_{4}$, or

$\qquad \qquad \qquad \qquad \qquad s_{3}-s_{4}=-1$ \ \ \ \ (2)

Equations (1) and (2) have the unique solution $s_{3}=3$, $s_{4}=4$, i.e., $%
G $ is the Borda correspondence.\medskip

Now suppose the Theorem holds for all $m<M$ and we have weights\medskip

\qquad \qquad $1<2\leq s_{3}\leq \cdots \leq s_{M}$\medskip \newline
such that for some $j>3$, it is the case that $s_{j}\neq j$. \ Let $G$ be
the corresponding scoring social choice correspondence. \ We show that there
is a profile $u$ such that $G(u)\neq G_{B}(u)$ and that $G$ must fail
balancedness. \ Consider the first $j$ such that $s_{j}\neq j$. \ If $j<M$,
just take a profile $v$ that works for $j$ alternatives (using the induction
hypothesis) and append $M-j$ additional alternatives to everyone's
bottom.\medskip

So we need only consider the case where the first $j$ such that $s_{j}\neq j$
is $j=M$. \ The weights are:\medskip

\qquad \qquad $1<2<3<\cdots <M-1\leq w$ where $w\neq M$.\medskip

Analysis here is for $w>M$. The same profile works for $w<M$.\medskip

We first show one stage. \ Look at the profile $u$ just above. Now insert
alternative $c$ just above $b$ for \#1 and just below $b$ for \#2. \ The
third individual has $c$ at the bottom:\medskip

\qquad \qquad \qquad \qquad 
\begin{tabular}{|l|l|l|}
\hline
$1$ & $2$ & $3$ \\ \hline
$a$ & \multicolumn{1}{|l|}{$b$} & $x$ \\ 
$x$ & \multicolumn{1}{|l|}{$c$} & $y$ \\ 
$y$ & \multicolumn{1}{|l|}{$a$} & $b$ \\ 
$c$ & \multicolumn{1}{|l|}{$y$} & $a$ \\ 
$b$ & $x$ & $c$ \\ \hline
\end{tabular}%
\medskip \newline
The Borda scores for $a$ and $x$ are the lowest, so $G_{B}(u)=\{a,x\}$. With
weights 1, 2, 3, 4, $w_{5}$, with $w_{5}>5$, the scores at $u$ are: $S(a)=8$%
, $S(b)=4+w_{5}$, $S(c)=6+w_{5}$, $S(x)=3+w_{5}$, and $S(y)=9$, so $%
G(u)=\{a\}$. \ After transposing $x$ and $y$ for individuals 1 and 2 to
create profile $v$, the scores become: $S(a)=8$, $S(b)=4+w_{5}$, $%
S(c)=6+w_{5}$, $S(x)=8$, and $S(y)=4+w_{5}$, so $G(v)=\{a,x\}$ and
balancedness is violated.\medskip

Eventually the profile looks like $u^{\ast }$:\medskip

\qquad \qquad \qquad \qquad 
\begin{tabular}{|l|l|l|}
\hline
$1$ & $2$ & $3$ \\ \hline
$a$ & \multicolumn{1}{|l|}{$b$} & $x$ \\ 
$x$ & $c_{M-4}$ & $y$ \\ 
$y$ & $\vdots $ & $b$ \\ 
$c_{1}$ & $c_{2}$ & $a$ \\ 
$c_{2}$ & $c_{1}$ & $c_{1}$ \\ 
$\vdots $ & $a$ & $c_{2}$ \\ 
$c_{M-4}$ & \multicolumn{1}{|l|}{$y$} & $\vdots $ \\ 
$b$ & $x$ & $c_{M-4}$ \\ \hline
\end{tabular}%
\medskip \newline
With weights 1, 2, ..., $M-1$, $w$, with $w>M$, the scores at $u^{\ast }$
are: $S(a)=M+3$, $S(b)=w+4$, $S(x)=w+3$, $S(y)=M+4$, and the rest have
higher scores. So $G(u^{\ast })=\{a\}$ given $w>M$. \ After transposing $x$
and $y$ for individuals 1 and 2 to create profile $v^{\ast }$, the scores
become: $S(a)=M+3$, $S(b)=w+4$, $S(x)=M+3$, $S(y)=w+4$, and the remaining
scores are unchanged, so $G(v^{\ast })=\{a,x\}$ and balancedness is
violated. \ \ \ \ $\square $\bigskip

We obtain from this result a new characterization of the Borda
correspondence. \ While Young (1974), Hansson and Sahlquist (1976), Coughlin
(1979/80), Nitzan and Rubinstein (1981), and Debord (1992) have
characterizations of Borda's rule, they work in a different context than
ours; for these authors, a rule has to work for a variable number of
individuals and they use variable population properties like that of
separability introduced by Smith (1973). Here we can get a characterization
of Borda's rule with a fixed set of individuals by appending balancedness to
a characterization of scoring rules for fixed populations [see Fishburn
(1973a, 1973b)].

\end{document}